\documentclass[a4paper,12pt]{article}

\usepackage{amsmath,amsfonts,amssymb,graphicx,color,array}

\newtheorem{newremark}{Remark}
\newenvironment{remark}
    {\begin{newremark}\rm}
    {\end{newremark}}
\newtheorem{theorem}{Theorem}

\newtheorem{newexample}{Example}
\newenvironment{example}
    {\begin{newexample}\it}
    {\end{newexample}}
\newtheorem{newlemma}{Lemma}
\newenvironment{lemma}
    {\begin{newlemma}\it}
    {\end{newlemma}}

\newcommand{\e}{{\rm e}}
\newcommand{\R}{\mathbb{R}}

\newcommand{\cF}{ {\cal F} }

\begin{document}

\title{Exponentially convergent method for integral nonlocal problem for the first order differential equation with unbounded coefficient in Banach space}

\author{V.B.~Vasylyk\thanks{Institute of Mathematics of NAS of Ukraine, 3 Tereshchenkivs'ka Str., Kyiv-4, 01601, Ukraine ({\tt vasylyk@imath.kiev.ua}).}  }

\date{\null}
\maketitle

\begin{abstract}
Problem for the first order differential equation with an unbounded operator coefficient in Banach space and integral nonlocal condition is considered. An exponentially convergent algorithm is proposed and justified for the numerical solution of this problem in assumption that an operator coefficient $A$ is strongly positive and some existence and uniqueness conditions are fulfilled. This algorithm  is based on the representations of operator functions by a Dunford-Cauchy integral along a hyperbola, enveloping the spectrum of $A$, and on the proper quadratures involving short sums of resolvents. The efficiency of the  proposed algorithms is demonstrated by several numerical examples.
\end{abstract}



\section{Introduction}

In this paper we consider the following nonlocal problem with integral condition:
\begin{equation}\label{pr-st}
  \begin{split}
     &\frac{du}{dt}+A u=0,  \quad t \in [0,T]\\
     &u(0)+\int_0^T w(s)u(s)ds =u_0,
  \end{split}
\end{equation}
where $w(s)\ge 0$ is a given function,  $u_0\in X.$ The operator $A$ with the domain $D(A)$ in a Banach space $X$ is assumed to be densely defined strongly positive (sectorial) operator, i.e. its spectrum $\Sigma(A)$ lies in a sector of the right half-plane with the vertex at the origin, while its resolvent decays inversely proportional to $|z|$ at the infinity (see estimate (\ref{estrez}) below).

 Inhomogeneous problem related to \eqref{pr-st} can be reduced to homogeneous one by change of function in the following way. If we have
\begin{equation}\label{pr-st-inh}
  \begin{split}
     &\frac{dv}{dt}+A v=f(t),  \quad t \in [0,T]\\
     &v(0)+\int_0^T w(s)v(s)ds =u_0,
  \end{split}
\end{equation}
with $f(t)$-- vector-valued function in the Banach space $X$ then putting $v(t)=u(t)+v_1(t)$, where
\[
v_1(t)=\int\limits_{0}^{t}\e^{-A(t-s)}f(s)ds,
\]
we obtain the following problem for $u(t)$
\[
  \begin{split}
     &\frac{du}{dt}+A u=0,  \quad t \in [0,T]\\
     &u(0)+\int_0^T w(s)u(s)ds =u_0-\Phi,
  \end{split}
\]
where 
\[
\Phi=\int\limits_{0}^{T}w(s)v_1(s)ds.
\]

Note that exponentially convergent numerical approximation for $v_1(t)$ was developed in \cite{gm5}, \cite{GMV-mon}. So, one can use this approximation to obtain $v_1(t)$ and then to find $\Phi$.

It should be noted that recently there were developed various exponentially convergent methods for problems with unbounded coefficients in Banach space \cite{fernandez-lubich-palencia-schaedle}, \cite{GMV-mon}, \cite{lopez-fernandez1}, \cite{sst1}, \cite{thomee1}, \cite{weid1}. These problems can be considered as metamodels of classical problems for partially differential equations such as parabolic elliptic and hyperbolic.

The aim of this paper is to construct an exponentially convergent approximation of a solution to problem \eqref{pr-st}. The paper is organized as follows. In Section \ref{stet} we discuss the existence and uniqueness of the solution as well as its representation through input data.  A numerical algorithm for the problem (\ref{pr-st}) is proposed and justified in section \ref{numalg}. The main result of this section is theorem \ref{oc-m} about the exponential convergence rate of the proposed discretization. The next section \ref{num-ex} is devoted to numerical examples which confirms theoretical results from the previous section.



\section{Existence and representation of the solution}\label{stet}

The solution of \eqref{pr-st} can be formally represented as follows
\begin{equation}\label{rr1}
u(t)= \e^{-At}u(0) .
\end{equation}
From the integral condition in \eqref{pr-st} and \eqref{rr1} we obtain
\[
u(0)+\int_0^T w(s)u(s)ds =u_0,
\]
\[
u(0)+\int_0^T w(s)\e^{-As}ds\, u(0)= u_0,
\]
Therefore, in the case when $\left[ I+ \int_0^T w(s)\e^{-As}ds\right] ^{-1}$ exists (sufficient conditions see below) we have
\[
u(0)=\left[ I+ \int_0^T w(s)\e^{-As}ds\right] ^{-1} u_0. 
\]
So,
\begin{equation}\label{reprsol}
u(t)=\e^{-At}\left[ I+ \int_0^T w(s)\e^{-As}ds\right] ^{-1} u_0 .
\end{equation}

Let the operator $A$ in \eqref{pr-st} be a densely defined strongly positive (sectorial) operator in a Banach space $X$ with the domain $D(A),$ i.e. its spectrum $\Sigma(A)$ is situated in a sector $\Sigma$
\begin{equation}\label{spectrA}
\Sigma = \left\{ z=\rho_0+r \e^{i\theta}:\quad r \in [0,\infty), \, \rho_0>0 \ \left|\theta\right|< \varphi < \frac{\pi}{2}
\right\}.
\end{equation}
Additionally,  the following estimate for the resolvent is valid
\begin{equation}\label{estrez}
\|R_A(z)\|=\left\|(zI-A)^{-1}\right\|\leq \frac{M}{1+\left|z\right|}
\end{equation}
outside the sector and on its boundary $\Gamma_\Sigma$. The numbers $\rho_0,$ $\varphi$ are called the spectral
characteristics of $A$.

We call a spectral hyperbola the curve $\Gamma_0:$
\begin{equation}\label{sp-hyp}
\Gamma_0=\{z(\zeta)=\rho_0 \cosh{\zeta}-i b_0 \sinh{\zeta}: \; \zeta \in
(-\infty, \infty ), \; b_0=\rho_0 \tan{\varphi}\}.
\end{equation}
It has a vertex at $(\rho_0,0)$ and asymptotes that are parallel to the rays of the spectral angle $\Sigma$.

A convenient representation of operator functions is the one through the Dunford-Cauchy integral (see e.g. \cite{Clem, krein}) where the integration path plays an important role. Using the Dunford-Cauchy integral representation and \eqref{reprsol} the solution to problem \eqref{pr-st} can be represented as
\begin{equation*}
     u(t)=\frac{1}{2\pi i}\int_{\Gamma_I} \e^{-zt}\left[ 1+ \int_0^T w(s)\e^{-zs}ds\right] ^{-1} R_A(z)u_0 dz=  
\end{equation*}
\[
=\frac{1}{2\pi i}\int_{\Gamma_I} F(z,A) R_A(z)u_0 dz ,
\]
if $F(z,A)$ is analytic function inside the integration hyperbola $\Gamma_I$ that envelopes $\Gamma_0.$ To obtain uniformly convergent and numerically stable algorithm we shall  modify this integral by changing the resolvent $R_A(z)$ to  $R_A^1(z)$ that doesn't change the value of integral when $u_0\in D(A^\alpha),$ $\alpha>0$ (for the details see \cite{GMV-mon}).
\[
R_A^1(z)= (zI-A)^{-1}-\frac{I}{z}.
\]
Therefore, one can obtain the following representation for the solution of the problem \eqref{pr-st}:
\begin{equation}\label{D-C-rep}
     u(t)=\frac{1}{2\pi i}\int_{\Gamma_I} F(z,A) R_A^1(z) u_0 dz  .
\end{equation}

 We choose the following hyperbola
\begin{equation}\label{int-hyp}
\Gamma_I=\{z(\zeta)=a_I \cosh{\zeta}-i b_I \sinh{\zeta}: \; \zeta \in (-\infty, \infty )\},
\end{equation}
for an integration contour that envelopes the spectrum of $A$, where the values of $a_I$, $b_I$ are to be defined later.  Using this hyperbola, we obtain  from \eqref{D-C-rep}
\begin{equation}\label{par-int}
 u(t)= \frac{1}{2 \pi i}\int_{-\infty}^{\infty}F(z(\zeta),A) R_A^1(\zeta)z^\prime(\zeta) u_0 d \zeta= \int_{-\infty}^{\infty}\cF(t,\zeta)d\zeta,
\end{equation}
with
\[
z^\prime(\zeta)=a_I \sinh{\zeta}-i b_I \cosh{\zeta}.
\]

The next step toward a numerical algorithm  is an approximation of \eqref{par-int} by the efficient quadrature formula. For this purpose we need to estimate the width of a strip around the real axis where the integrand in \eqref{par-int} permits analytical extension (with respect to $\zeta$).  The integration hyperbola $\Gamma_I$ will be translated into the parametric set of hyperbolas with respect to $\nu$ after changing $\zeta$ to $\zeta+i \nu$
\begin{equation*}\label{par-hyp}
\begin{split}
\Gamma(\nu)&=\{z(\zeta,\nu)=a_I \cosh{(\zeta+i \nu)}-ib_I \sinh{(\zeta+ i \nu)}: \; \zeta\in(-\infty,\infty)\} \\
&=\{z(\zeta,\nu)=a(\nu) \cosh{\zeta}-ib(\nu) \sinh{\zeta}: \; \zeta\in(-\infty,\infty)\},
\end{split}
\end{equation*}
with
\begin{equation*}\label{na12}
\begin{split}
&a(\nu)= a_I \cos{\nu}+b_I
\sin{\nu}=\sqrt{a_I^2+b_I^2}\sin{(\nu+\phi/2)}, \\
& b(\nu)= b_I \cos{\nu}-a_I
\sin{\nu}=\sqrt{a_I^2+b_I^2}\cos{(\nu+\phi/2)},\\
& \cos{\frac{\phi}{2}}=\frac{b_I}{\sqrt{a_I^2+b_I^2}}, \;
\sin{\frac{\phi}{2}}=\frac{a_I}{\sqrt{a_I^2+b_I^2}} \; .
\end{split}
\end{equation*}

The analyticity of the integrand in the strip
\begin{equation*}\label{na16}
\begin{split}
&D_{d_1}=\{(\zeta, \nu):\zeta \in (-\infty, \infty), |\nu|<d_1/2\},
\end{split}
\end{equation*}
with some $d_1$ could be violated if the resolvent or the part related to the nonlocal condition become unbounded. To avoid this we have to choose $d_1$ in a way such that for $\nu \in (-d_1/2,d_1/2)$ the hyperbola $\Gamma(\nu)$ remains in the right half-plane of the complex plane. For $\nu=-d_1/2$ the corresponding hyperbola is going through the point $(\rho_1,0)$, for some $0 \le \rho_1 <\rho_0$. For $\nu=d_1/2$ it coincides with the spectral hyperbola and therefore for all $\nu \in (-d_1/2,d_1/2)$ the set {$\Gamma(\nu)$} does not intersect the spectral sector. This fact justifies the choice the hyperbola $\Gamma(0)=\Gamma_I$ as the integration path.

Such requirements for $\Gamma(\nu)$ imply the following system of equations
\begin{equation*}
  \begin{cases}
    a_I \cos{(d_1/2)}+b_I \sin{(d_1/2)}=\rho_0 ,   \\
   b_I \cos{(d_1/2)}-a_I \sin{(d_1/2)}=b_0=\rho_0\tan{\varphi} ,  \\
   a_I \cos{(-d_1/2)}+b_I \sin{(-d_1/2)}=\rho_1 ,
  \end{cases}
\end{equation*}
it leads us to the next system
\begin{equation*}
  \begin{cases}
   a_I=\rho_0 \cos{(d_1/2)} -b_0 \sin{(d_1/2)},  \\
   b_I=\rho_0 \sin{(d_1/2)} +b_0 \cos{(d_1/2)}, \\
   2 a_I \cos{(d_1/2)}=\rho_0+ \rho_1.
  \end{cases}
\end{equation*}
Eliminating $a_I$ from the first and the third equations we obtain
\[
\rho_0 \cos{d_1}- b_0 \sin{d_1}=\rho_1,
\]
\[
\cos(d_1+\varphi)=\frac{\rho_1}{\sqrt{\rho_0^2+b_0^2}},
\]
i.e.
\begin{equation}\label{shyr-sm}
d_1=\arccos{\left(\frac{\rho_1}{\sqrt{\rho_0^2+b_0^2}}\right)} -\varphi ,
\end{equation}
with $\cos{\varphi}=\frac{\rho_0}{\sqrt{\rho_0^2 +b_0^2}},$ $\sin{\varphi}=\frac{b_0}{\sqrt{\rho_0^2+b_0^2}}.$ Thus,
for $a_I$, $b_I$ we receive
\begin{equation}\label{na15}
\begin{split}
a_I&=\sqrt{\rho_0^2+b_0^2} \cos{\left(\frac{d_1}{2} +\varphi\right)} \\
& = \rho_0 \frac{\cos{\left(\frac{d_1}{2} +\varphi\right)}}{\cos{\varphi}} =\rho_0 \frac{\cos{\left(\arccos\left(\frac{\rho_1}{\sqrt{\rho_0^2
+b_0^2}}\right)/2  +\varphi/2 \right)}}{\cos{\varphi}}, \\
b_I&=\sqrt{\rho_0^2+b_0^2} \sin{\left(\frac{d_1}{2} +\varphi\right)} \\
&=\rho_0 \frac{\cos{\left(\frac{d_1}{2} +\varphi\right)}}{\cos{\varphi}} =\rho_0 \frac{\cos{\left(\arccos\left(\frac{\rho_1}{\sqrt{\rho_0^2
+b_0^2}}\right)/2  +\varphi/2 \right)}}{\cos{\varphi}}.
\end{split}
\end{equation}
For $a_I$ and $b_I$ defined as above the resolvent of the operator $A$ is analytic in the strip $D_{d_1}$ with respect to $w=\zeta+i \nu$  for any $t \ge 0$. Note, that for $\rho_1=0$ we have $d_1=\pi/2-\varphi$ as in \cite{gm5}.

Taking into account (\ref{na15}) we can similarly write equations for $a(\nu),$ $b(\nu)$ on the whole interval
$-\frac{d_1}{2}\le\nu \le \frac{d_1}{2}$
\[
\begin{split}
a(\nu) &= a_I \cos{\nu}+ b_I\sin{\nu} =\sqrt{\rho_0^2+b_0^2}\cos{\left(\frac{d_1}{2} +\varphi\right)} \cos(\nu) \\
 &+\sqrt{\rho_0^2+b_0^2}\sin{\left(\frac{d_1}{2}+\varphi\right)} \sin(\nu) =\sqrt{\rho_0^2+b_0^2}\cos{\left(\frac{d_1}{2}+\varphi-\nu\right) } , \\
 b(\nu) &= b_I \cos{\nu}-a_I \sin{\nu} =\sqrt{\rho_0^2+b_0^2}\sin{\left(\frac{d_1}{2} +\varphi\right)} \cos(\nu) \\
 &-\sqrt{\rho_0^2+b_0^2}\cos{\left(\frac{d_1}{2}+\varphi\right)} \sin(\nu) =\sqrt{\rho_0^2+b_0^2}\sin{\left(\frac{d_1}{2}+\varphi -\nu\right)} ,
\end{split}
\]
\[
 \rho_1\le a(\nu)\le \rho_0, \quad b_0\le b(\nu)\le \sqrt{b_0^2+\rho_0^2 -\rho_1^2},
\]
with $d_1,$ defined by (\ref{shyr-sm}).

Now, let us establish condition on $w(s)$, when expression
\[
\left[ 1+ \int_0^T w(s)\e^{-zs}ds\right]
\]
related to nonlocal condition dose not become zero inside the integration hyperbola $\Gamma_I$. 
\[
\left| 1+ \int_0^T w(s)\e^{-z(\zeta)s}ds\right| \ge 1- \left| \int_0^T w(s)\e^{-z(\zeta)s}ds\right| \ge
\]
\[
\ge 1- \|w(s)\|_{C[0,T]} \int_0^T \e^{-sa_I\cosh(\zeta)}ds= 1-\frac{\|w(s)\|_{C[0,T]}}{a_I \cosh(\zeta) }\left( 1-\e^{-Ta_I\cosh(\zeta)}\right) \ge
\]
\[
\ge 1-\frac{\|w(s)\|_{C[0,T]}}{a_I }.
\]
Therefore, we have
\[
\left| 1+ \int_0^T w(s)\e^{-z(\zeta)s}ds\right|^{-1}\le \frac{1}{1-\frac{\|w(s)\|_{C[0,T]}}{a_I }} \le C_1,
\]
in the case when
\begin{equation}\label{umzb}
\|w(s)\|_{C[0,T]}<a_I,
\end{equation}
where $a_I$ is defined in \eqref{na15}.

So, we can summarize all of the above in the following lemma.
\begin{lemma}
Let $A$ be a densely defined strongly positive operator. If the condition \eqref{umzb} is valid then there exists a unique solution to problem \eqref{pr-st} that can be represented by \eqref{D-C-rep}.
\end{lemma}

More rough estimate than \eqref{umzb} is
\begin{equation}\label{umzb2}
\|w(s)\|_{C[0,T]}\le \frac{1}{T},
\end{equation}
that can be easily obtained from estimate
\[
\left| \int_0^T w(s)\e^{-z(\zeta)s}ds\right|\le \|w(s)\|_{C[0,T]}T.
\]

Further, let us establish conditions for the existence of the solution to \eqref{pr-st} in the case when the  operator $A$ is self-adjoint positive definite. To achieve that we have to choose $d_1$ in a way that for $\nu \in (-d_1/2,d_1/2)$ the hyperbola $\Gamma(\nu)$ remains in the right half-plane of complex plane. For $\nu=-d_1/2$ the corresponding hyperbola turns into the line parallel to the imaginary axis. For $\nu=d_1/2$ it coincides with the ray that lies on the real axis having a vertex at $\rho_0$. These requirements imply the following system of equations
\begin{equation*}
  \begin{cases}
   a_I \cos{(d_1/2)}+b_I \sin{(d_1/2)}=\rho_0 ,   \\
   b_I \cos{(d_1/2)}-a_I \sin{(d_1/2)}=0 ,  \\
   a_I \cos{(-d_1/2)}+b_I \sin{(-d_1/2)}=0 ,
  \end{cases}
\end{equation*}
which has the solution
\[
a_I=b_I=\frac{\rho_0}{\sqrt{2}},
\]
\[
d_1=\frac{\pi}{2}
\]
The condition \eqref{umzb} then  becomes
\begin{equation}\label{umM2}
    \|w(s)\|_{C[0,T]} < \frac{\rho_0}{\sqrt{2}},
\end{equation}
that is sufficient condition of existence solution to \eqref{pr-st} in the case of self-adjoint positive operator $A$.

\section{Numerical algorithm}\label{numalg}

First of all we approximate integral in denominator using exponentially convergent quadrature. For this aim we use Gauss quadrature rule. So,
\begin{equation}\label{Gquad}
I=\int_0^T w(s)\e^{-z(\zeta)s}ds\approx \sum\limits_{j=0}^{n}\frac{T}{2}\omega_j w(\xi_j)\e^{-z(\zeta)\xi_j}=I_n,
\end{equation}
\[
\xi_j=\frac{T}{2}(\theta_j+1),
\]
where $\theta_j$-- are a set of $n+1$  roots of the Legendre polynomial $P_{n+1}(x)$ and $\omega_j$-- a are a set of weights related to the Gauss quadrature rule. Note that $\theta_j$ and $\omega_j$ can be precomputed using fast algorithms (see \cite{TrefApp}).

Therefore we obtain from \eqref{par-int}
\begin{equation}\label{intnab}
u(t)\approx u_n(t)= \frac{1}{2 \pi i}\int_{-\infty}^{\infty}F_n(z(\zeta),A) R_A^1(\zeta)z^\prime(\zeta) u_0 d \zeta=\int_{-\infty}^{\infty}\cF_n(t,\zeta) d \zeta,
\end{equation}
where
\[
F_n(z(\zeta),A)=\e^{-z(\zeta)t}\left[ 1+ I_n \right] ^{-1}
\]

For the error estimate we have
\[
\left| \frac{1}{1+I}-\frac{1}{1+I_n} \right| = \left| \frac{I_n-I}{(1+I)(1+I_n)} \right|. 
\]
Due to \eqref{umzb} or \eqref{umzb2} we have
\[
 \frac{1}{\left| 1+I\right|} \le C. 
\]
For the second multiplier we obtain
\[
\frac{1}{\left| 1+I_n\right|} \le \frac{1}{1-\left| \frac{T}{2}\sum\limits_{j=0}^{n}\omega_j w(\xi_j)\e^{-z(\zeta)\xi_j}\right| } \le \frac{1}{1-\frac{T\|w(s)\|_{C[0,T]}}{2}  \sum\limits_{j=0}^{n}\omega_j \e^{-a_I\cosh(\zeta)\xi_j} }\le 
\]
\begin{equation}\label{estIntn}
\le \frac{1}{1-T \|w(s)\|_{C[0,T]} }\le c=const,
\end{equation}
in the case when \eqref{umzb2} is valid. Therefore we have
\[
\left| \frac{1}{1+I}-\frac{1}{1+I_n} \right| \le c \left| I_n-I \right|. 
\]
Exponential function $\e^{-zs}$ is analytical on $s$ in all complex plane. So, smoothness of integrand in $I$ is defined by $w(s)$. Using theorem 19.3 from \cite{TrefApp} we have that if $w(\frac{T}{2}(s+1))$ is analytic in $[-1,1]$, analytically continuable to the open Bernstein ellipse where $\left| w(\frac{T}{2}(s+1))\e^{-z\frac{T}{2}(s+1)}\right| \le M$ then
\begin{equation}\label{intAn}
 \left| I-I_n\right| \le \frac{144M\rho^{-2n}}{35(\rho^2-1)}, \qquad n\ge 2. 
\end{equation}
If $w(s)$ and its derivatives through $w^{(\nu-1)}$ are absolutely continuous and $w^{(\nu)}$ is of bounded variation $V$ then
\begin{equation}\label{intdiff}
 \left| I-I_n\right| \le \frac{32V}{15\pi \nu (n-2\nu-1)^{2\nu+1}}, \qquad n> 2\nu+1 . 
\end{equation}

Supposing $u_0 \in D(A^{\alpha}),$ $0<\alpha<1$ it was shown in \cite{GMV-mon} that
\begin{equation}\label{dopest}
\begin{split}
  \left\| \e^{-z(\zeta)t}z^\prime(\zeta) R_A^1(\zeta) u_0 \right\| &
  \le (1+M)K\frac{b_I}{a_I}\left(\frac{2}{a_I}\right)^{\alpha} \e^{-a_I t \cosh{\zeta}-\alpha |\zeta|} \|A^{\alpha}u_0\|, \\
& \zeta \in \R, \quad t\ge 0,
\end{split}
\end{equation}
where $K$ is a constant that depends on $\alpha$, $M$ is a constant from resolvent estimate \eqref{estrez}.

The part responsible for the nonlocal condition in \eqref{intnab} is estimated by \eqref{estIntn}. Thus, we obtain the following estimate for $\cF(t,\xi)$:
\begin{equation}\label{oc-pid}
\begin{split}
    \|\cF_n(t,\zeta) \| \le C(\varphi,\alpha) \e^{-a_I t \cosh{\zeta}-\alpha |\zeta|} \|A^{\alpha}u_0\|, \\
 C(\varphi,\alpha)=\frac{(1+M)Kc b_I}{2\pi a_I}\left(\frac{2}{a_I}\right)^{\alpha},\quad \zeta \in \R, \quad t\ge 0.
\end{split}
\end{equation}

We approximate  integral (\ref{intnab}) by the following Sinc-quadrature formula \cite{GMV-mon, stenger}:
\begin{equation}\label{h-nab}
u_{n,N}(t)=h\sum_{k=-N}^{N}\cF_n(t,z(kh)),
\end{equation}
with the error
\[
\|\eta_N(\cF_n,h)\|=\|u_n(t)-u_{n,N}(t)\|
\]
\[
\le \left\|u_n(t)- h \sum_{k=-\infty}^{\infty}\cF_n(t,z(kh))\right\| +\left\|h\sum_{|k|>N}\cF_n(t,z(kh))\right\|
\]
\[
\le \frac{1}{4 \pi}\frac{\e^{-\pi d/h}}{ \sinh{(\pi d/h)}}\|\cF_n\|_{{\bf H}^1(D_{d})}
\]
\[
+C(\varphi,\alpha)h \|A^{\alpha}u_0\| \sum_{k=N+1}^{\infty} \e^{-a_I t \cosh{kh}-\alpha kh}.
\]

Here $ {\bf H}^1(D_d)$ is a space introduced similarly to \cite{stenger} in  \cite{GMV-mon} of all vector-valued functions $\cal{F}$ analytic in the strip $D_d$ . Due to \cite{GMV-mon}

\begin{equation}\label{na18-1-1}
\begin{split}
& \|\e^{-z(\cdot)t}z^\prime(\cdot) R_A^1(\cdot) u_0\|_{{\bf H}^1(D_d)}\le  \|A^{\alpha}u_0\|[C_{-}(\varphi,\alpha, \delta)\\
 &+C_{+}(\varphi,\alpha, \delta)]
 \int_{-\infty}^\infty \e^{-\alpha |\xi|} d
 \xi=C(\varphi,\alpha,\delta)\|A^{\alpha}u_0\|
 \end{split}
\end{equation}
with
\begin{equation}\label{na19-1-1}
\begin{split}
& C(\varphi,\alpha,\delta)=
\frac{2}{\alpha}[C_{+}(\varphi,\alpha,\delta)+C_{-}(\varphi,\alpha,\delta)], \\
&C_{\pm}(\varphi,\alpha,\delta)=(1+M)K
\tan{\left(\frac{\pi}{4}+\frac{\varphi}{2}\pm \frac{d}{2} \right)}
\left(\frac{2\cos{\varphi}}{a_0\cos{\left(\frac{\pi}{4}+\frac{\varphi}{2}\pm
\frac{d}{2} \right)}}\right)^\alpha,
\end{split}
\end{equation}
\[
d=d_1-\delta,
\]
for an arbitrary small positive $\delta$. 

Obviously that in the case of  \eqref{umzb2} the part responsible for the nonlocal condition  is bounded in $D_{d}$. Therefore 
we obtain
\[
\|\cF_n(t,\zeta)\|_{{\bf H}^1(D_d)} \le C(\varphi,\alpha,\delta)\|A^{\alpha}u_0\|.
\]

So, we have for the error $\eta_N(\cF_n,h)$
\begin{equation}\label{err1}
\|\eta_N(\cF_n,h)\| \le \frac{c \|A^{\alpha}u_0\|}{\alpha}\left\{\frac{\e^{-\pi d/h}}{
\sinh{(\pi d/h)}}+\e^{-a_I t \cosh{((N+1)h)}-\alpha (N+1) h}\right\}
\end{equation}
where the constant $c$ does not depend on $h$, $N$, $t$.

Equalizing the both exponentials gives 
\begin{equation*}
 \frac{ \pi d}{h}=\alpha(N+1)h,
\end{equation*}
\begin{equation}\label{na23}
 h=\sqrt{\frac{ \pi d_1}{\alpha(N+1)}},
\end{equation}
this leads us to the following error estimate
\begin{equation}\label{na24}
\|\eta_N(\cF_n,h)\|\le \frac{c}{\alpha}\text{exp}{\left(-\sqrt{\pi d \alpha(N+1)}\right)}\|A^{\alpha}u_0\|
\end{equation}
In the case $t >1$ the first summand in the argument of $\e^{-a_I t \cosh{((N+1)h)}-\alpha (N+1) h}$ in \eqref{err1} contributes mainly to the error. Setting for such case $h=c_1 \ln{N}/N$ with some positive constant $c_1$ we obtain for a fixed $t$ the following estimate:
\begin{equation}\label{na240}
\|\eta_N(\cF_n,h)\|\le c\left[\e^{- \pi d_1 N/(c_1 \ln{N})} +\e^{-c_1 t a_I N/2-c_1 \alpha \ln{N} }\right] \|A^{\alpha}u_0\|.
\end{equation}
Thus, we have proven the following theorem.
\begin{theorem}\label{oc-m}
Let $A$ be a densely defined strongly positive operator, $u_0 \in D(A^{\alpha})$, $\alpha \in (0,1)$ and condition \eqref{umzb2} is valid. Then Sinc-quadrature \eqref{h-nab} represents an approximate to $u_n(t)$. It provides the convergence of exponential order uniformly with respect to $t \ge 0$ presented by the estimate \eqref{na24} for the step size $h$ defined in \eqref{na23}. The approximation has the convergence rate \eqref{na240} for the case $t>0$ and $h=c_1 \ln{N}/N$.
\end{theorem}

\begin{remark}
 The integration curve $\Gamma_I$ is symmetric with respect to the real axis. Therefore $z(-kh)=\overline{z(kh)}$ and $z'(-kh)=-\overline{z'(kh)}$. Approximation \eqref{h-nab} can be rewritten in the form
 \[
 u_{n,N}(x)=\frac{h}{2 \pi i}\cF_n(t,z(0)) + Re \left[ \sum_{k=1}^{N}h \frac{\cF_n(t,z(kh))}{\pi i}\right],
 \]
 which reduce the number of  resolvent calculations by factor of two.
\end{remark}

Further we consider the full error estimate.
\[
\varepsilon_1=\|u(t)-u_n(t) \|=\left\| \int_{-\infty}^{\infty}\left[ \cF(t,\zeta)-\cF_n(t,\zeta)\right] d\zeta \right\| \le
\]
\[
\frac{1}{2\pi}\int_{-\infty}^{\infty}\left| \e^{-z(\zeta)t}\zeta^\prime(\zeta)\right| \left| \frac{1}{1+I}-\frac{1}{1+I_n} \right| \left\| R_A^1(\zeta) u_0\right\|  d\zeta,
\]
then using \eqref{dopest} we obtain
\[
\varepsilon_1= \frac{(1+M)Kb_Ic}{2\pi a_I}\left( \frac{2}{a_I}\right)^\alpha \left\| A^\alpha u_0\right\| \left| I-I_n\right| \int_{-\infty}^{\infty} \e^{-a_I t\cosh\zeta -\alpha|\zeta|} d\zeta  \le
\]
\[
\le \frac{(1+M)Kb_Ic}{\pi a_I}\left( \frac{2}{a_I}\right)^\alpha \left\| A^\alpha u_0\right\| \left| I-I_n\right| \int_{0}^{\infty} \e^{-\alpha|\zeta|} d\zeta=
\]
\[
= \frac{(1+M)Kb_Ic}{\pi a_I\alpha}\left( \frac{2}{a_I}\right)^\alpha \left\| A^\alpha u_0\right\| \left| I-I_n\right|= C\left\| A^\alpha u_0\right\| \left| I-I_n\right|.
\]

For the full error we have
\[
\|u(t)-u_{n,N}(t) \|\le \varepsilon_1+\left\| \eta_N(\cF_n,h) \right\| .
\]

Therefore we can formulate the main theorem.

\begin{theorem}\label{mm}
Let the conditions of theorem \ref{oc-m} are valid. Then \eqref{h-nab} represents an approximate to $u(t)$. It provides the convergence of exponential order in the case when $w(t)$ is analytically continuable to the Bernstein ellipse.
\end{theorem}


\section{Numerical examples}\label{num-ex}
\begin{example}\label{ex1}
Let us consider the problem \eqref{pr-st} with the operator $A$ defined  by
\begin{equation}\label{oper-def}
\begin{split}
&D(A)=\{v(x) \in H^2(0,1): \; v(0)=v(1)=0\}, \; \\
&Av=-v^{\prime \prime}(x) \; \forall v \in D(A),
\end{split}
\end{equation}
that generates a homogeneous parabolic equation with boundary conditions
\[
\begin{split}
 &\frac{\partial u(x,t)}{\partial t} -\frac{\partial ^2u(x,t)}{\partial x^2}=0,\\
 &u(0,t)=u(1,t)=0.
\end{split}
\]
Let us supplement this problem with a nonlocal integral condition
\[
u(x,0)+\int_{0}^{\frac{\pi}{2}}\cos (s) u(x,s)ds= \frac{\pi^4+\pi^2+\e^{-\pi^3/2}}{\pi^4+1}\sin(\pi x).
\]

In this case the exact solution to the problem is $u(x,t)=\e^{-\pi^2t}\sin(\pi x).$ We have  performed calculations using Maple. The errors are presented in Tables \ref{tab1}, \ref{tab2} for different $n$-- number of quadrature points \eqref{Gquad} and $N$-- number of Sinc-points \eqref{h-nab}. The table clearly exhibits an exponential decay of error according to the theoretical estimate.
\begin{table}[ptbh]
\begin{center}
  \begin{tabular}{|c||c|c|c|}
    \hline
    & \multicolumn{3}{c|}{n}\\ 
   \cline{2-4} N & 4 & 8 & 16 \\   \hline
    $4$  &  $0.001074121004$      &   & \\
    $8$  &  $0.000004530997940$   & $0.00000418248071$           & \\
    $16$ &  $2.394152400*10^{-7}$ & $7.3086845013760*10^{-10}$   & $7.159165797001*10^{-10}$\\
    $32$ &  $2.387505905*10^{-7}$ & $8.2307398421915*10^{-12}$   & $2.609087146562*10^{-13}$\\
    $64$ &                        & $7.9836951021369*10^{-12}$   & $2.917976861643*10^{-18}$\\
    $128$&                        &                              & $1.622297889726*10^{-24}$\\
    $256$&                        &                              & $1.873772451287*10^{-24}$\\
                \hline
  \end{tabular}
\end{center}
\caption{The error for $x=0.5, \; t=1$}\label{tab1}
\end{table}
\begin{table}[ptbh]
\begin{center}
  \begin{tabular}{|c||c|c|}
    \hline
    & \multicolumn{2}{c|}{n}\\ 
   \cline{2-3} N & 32 & 64  \\   \hline
    $128$  &  $2.559484448336975*10^{-25}$      &   \\
    $256$  &  $2.398463850885652*10^{-35}$      & $2.3984646885635428*10^{-35}$          \\
    $512$ &  $1.564339690250043*10^{-49}$       & $1.5716982365989911*10^{-49}$  \\
    $1024$ &                                    & $8.2307398421915*10^{-69}$  \\
                 \hline
  \end{tabular}
\end{center}
\caption{The error for $x=0.5, \; t=1$}\label{tab2}
\end{table}
\end{example}

\begin{example}
In this example we consider the problem as in example \ref{ex1} with a nonlocal integral condition
\[
u(x,0)+\int_{0}^{\frac{\pi}{2}}\cos (s^2) u(x,s)ds= (1-x)x^2.
\]

The results of calculation are presented in table \ref{tab3} for different $n$-- number of quadrature points \eqref{Gquad} and $N$-- number of Sinc-points \eqref{h-nab}. The table exhibits an exponential decay of error according to the theoretical estimate.
\begin{table}[ptbh]
\begin{center}
  \begin{tabular}{|c||c|c|c|}
    \hline
    n/N  & $u(x,t)$ \\   \hline
    $n=4$,   $N=32$   &  $0.5979651691*10^{-4}$           \\
    $n=8$,   $N=64$   &  $0.595184687264196427200402957709*10^{-4}$          \\
    $n=16$,  $N=128$  &  $0.595184553823189342113182135931*10^{-4}$        \\
    $n=32$,  $N=256$  &  $0.595184553823189342143429937050*10^{-4}$        \\
    $n=64$,  $N=512$  &  $0.595184553823189342143429937049*10^{-4}$    \\
    $n=128$, $N=1024$ &  $0.595184553823189342143429937049*10^{-4}$    \\
    \hline
  \end{tabular}
\end{center}
\caption{Solution to the problem for $x=0.4, \; t=1$}\label{tab3}
\end{table}

\end{example}


\begin{thebibliography}{10}

\bibitem{Clem}
Ph. Clement, H.J.A.M. Heijmans, S.~Angenent, C.J. van Duijn, and B.~de~Pagter.
\newblock {\em One-parameter semigroups.}
\newblock CWI Monographs, 5. North-Holland Publishing Co., Amsterdam, 1987.

\bibitem{fernandez-lubich-palencia-schaedle}
M.L. Fernandez, Ch. Lubich, C.~Palencia, and A.~Sch{\"a}dle.
\newblock Fast {R}unge-{K}utta approximation of inhomogeneous parabolic
  equations.
\newblock {\em Numerische Mathematik}, pages 1--17, 2005.

\bibitem{GMV-mon}
I.~Gavrilyuk, V.~Makarov, and V.~Vasylyk.
\newblock {\em Exponentially convergent algorithms for abstract differential
  equations}.
\newblock Frontiers in Mathematics. Birkh\"auser/Springer Basel AG, Basel,
  2011.

\bibitem{gm5}
I.P. Gavrilyuk and V.L. Makarov.
\newblock Exponentially convergent algorithms for the operator exponential with
  applications to inhomogeneous problems in {B}anach spaces.
\newblock {\em {S}{I}{A}{M} {J}ournal on Numerical Analysis}, 43(5):2144--2171,
  2005.

\bibitem{krein}
S.G. Krein.
\newblock {\em Linear Differential Operators in Banach Spaces}.
\newblock Amer. Math. Soc., New York, 1971.

\bibitem{lopez-fernandez1}
M.~L{\'o}pez-{F}ern{\'a}ndez, C.~Palencia, and A.~Sch{\"a}dle.
\newblock A spectral order method for inverting sectorial laplace transforms.
\newblock {\em SIAM J. Numer. Anal.}, 44:1332--1350, 2006.

\bibitem{sst1}
D.~Sheen, I.~H. Sloan, and V.~Thom{\'e}e.
\newblock A parallel method for time-discretization of parabolic equations
  based on laplace transformation and quadrature.
\newblock {\em IMA Journal of Numerical Analysis}, 23:269--299, 2003.

\bibitem{stenger}
F.~Stenger.
\newblock {\em Numerical methods based on Sinc and analytic functions}.
\newblock Springer Verlag, New York, Berlin, Heidelberg, 1993.

\bibitem{thomee1}
V.~Thom{\'e}e.
\newblock A high order parallel method for time discretization of parabolic
  type equations based on {L}aplace transformation and quadrature.
\newblock {\em Int. J. Numer. Anal. Model.}, 2:121--139, 2005.

\bibitem{TrefApp}
Lloyd~N. Trefethen.
\newblock {\em Approximation Theory and Approximation Practice}.
\newblock SIAM, 2012.

\bibitem{weid1}
J.~A.~C. Weideman.
\newblock Optimizing talbot's contours for the inversion of the laplace
  transform.
\newblock {\em SIAM J. Numer. Anal.}, 44(6):2342--2362, 2006.

\end{thebibliography}


\end{document}